\documentclass{article}

\usepackage{amsmath}
\usepackage{amssymb}

\newtheorem{thm}{Theorem}

\newtheorem{cor}{Corollary}

\title{Numbers with three close factorizations}
\author{Tsz Ho Chan}
\date{}

\begin{document}
\maketitle

\begin{abstract}
In this paper, we study numbers $n$ that can be factored in three different ways as $n = A_1 B_1 = A_2 B_2 = A_3 B_3$ where $A_1$, $A_2$, $A_3$ are close to each other and $B_1$, $B_2$, $B_3$ are close to each other.
\end{abstract}

\section{Introduction and main results}
Suppose a positive integer $n$ can be factored as $n = A B$. Is it possible to have another factorization $n = (A + a)(B - b)$ with small positive integers $a$ and $b$? The answer is yes. For this to hold, we require $A B = (A + a) (B - b)$ or $a b = a B - b A$. Let $d = (a, b)$ and $a = d a'$, $b = d b'$. Then $d a' b' = a' B - b' A$. This implies that $a' | A$ and $b' | B$ since $a'$ and $b'$ are relatively prime. Let $A = a' A'$ and $B = b' B'$. Then $d = B' - A'$. Therefore any such number would have the form $n = (a' A')(b' (A' + d)) = (a' (A' + d))(b' A')$ where $a'$, $b'$, $d$ are small and $A'$ can be arbitrarily large. What if we ask for two such extra factorizations, i.e. $n = A B = (A + a_1) (B - b_1) = (A + a_2) (B - b_2)$? It turns out that $A$ and $B$ cannot be arbitrarily large. We have
\begin{thm} \label{thm1}
Given $C \ge 2$. If $n = A B = (A + a_1)(B - b_1) = (A + a_2)(B - b_2)$ with $1 \le a_1 < a_2 \le C$ and $1 \le b_1 < b_2 \le C$, then $A, B < C^3$.
\end{thm}
One can ask if the upper bound $C^3$ is sharp and we have the following.
\begin{thm} \label{thm2}
Given $C \ge 10$. If $n = A B = (A + a_1)(B - b_1) = (A + a_2)(B - b_2)$ with $1 \le a_1 < a_2 \le C$ and $1 \le b_1 < b_2 \le C$, then $A, B \le \frac{1}{4} C^3 - \frac{1}{2} C^2 + \frac{1}{4} C$.
\end{thm}
Here the upper bound is best possible. For example, based on the proof of Theorem \ref{thm2} (with $a_1 = N$, $a_2 = 2N+1$, $b_1 = N+1$, $b_2 = 2N+3$), consider
\begin{equation} \label{eq1}
n = N (N+1)^2 (N+2) (2N+1) (2N+3).
\end{equation}
Firstly,
\[
n = [(2N+1) N (N+2)] \cdot [(2N+3) (N+1)^2] =: A \cdot B.
\]
Secondly,
\[
n = [(2N+3)(N+1)N] \cdot [(2N+1)(N+1)(N+2)] = [(2N+1) N (N+2) + N] \cdot [(2N+3) (N+1)^2 - (N+1)].
\]
Thirdly,
\[
n = [(2N+1) (N+1)^2] \cdot [(2N+3) N (N+2)] = [(2N+1) N (N+2) + (2N+1)] \cdot [(2N+3) (N+1)^2 - (2N+3)].
\]
With $C = 2N+3$, this shows that there are infinitely many integers $n$ satisfying the hypothesis in Theorem \ref{thm1} or \ref{thm2} with $\max(A, B) \ge \frac{1}{4} C^3 - \frac{1}{2} C^2 + \frac{1}{4} C$.

\bigskip

In a different perspective, we can think of finding three close factorizations of $n$ as finding three close lattice points on the hyperbola $x y = n$. Granville and Jim\'{e}nez-Urroz [\ref{GJ}] gave a lower bound for an arc of the hyperbola containing $k$ integer lattice points. Suppose $(x_1, y_1)$, $(x_2, y_2)$ and $(x_3, y_3)$ are three integer lattice points on $x y = n$ with $x_1 < x_2 < x_3$. They showed that $x_3 - x_1 \ge 2^{2/3} x_1 / n^{1/3}$. From Theorem \ref{thm2}, we have
\begin{cor} \label{cor1}
Suppose $(x_1, y_1)$, $(x_2, y_2)$ and $(x_3, y_3)$ are three integer lattice points on $x y = n$ with $x_1 < x_2 < x_3$. Then $\max(x_3 - x_1, y_1 - y_3) > 2^{2/3} n^{1/6} + 0.5$.

Moreover there exist infinitely many integer $n$ such that the hyperbola $x y = n$ contains three integer lattice points $(x_1, y_1)$, $(x_2, y_2)$ and $(x_3, y_3)$ with $x_1 < x_2 < x_3$ satisfying $\max(x_3 - x_1, y_1 - y_3) < 2^{2/3} n^{1/6} + 1.2$.
\end{cor}
This shows that the constant $2^{2/3}$ in all of the above is sharp. It also shows that the lower bound in the first half of Corollary \ref{cor1} is best possible (apart from a constant) which vanquishes any hope to improve the lower bound to $\gg n^{1/4}$ as claimed in [\ref{GJ}].

\section{Proof of Theorem \ref{thm1}}
Suppose $A B = (A + a_1)(B - b_1) = (A + a_2)(B -  b_2)$ with $1 \le a_1 < a_2 \le C$ and $1 \le b_1 < b_2 \le C$. Then $a_1 B - b_1 A = a_1 b_1$. Dividing by $b_1 B$, we have $\frac{a_1}{b_1} - \frac{A}{B} = \frac{a_1}{B}$. Similarly $\frac{a_2}{b_2} - \frac{A}{B} = \frac{a_2}{B}$. Subtracting the two equations, we have $\frac{a_2}{b_2} - \frac{a_1}{b_1} = \frac{a_2 - a_1}{B} > 0$. Hence $\frac{1}{C^2} \le \frac{1}{b_1 b_2} \le \frac{a_2 b_1 - a_1 b_2}{b_1 b_2} = \frac{a_2 - a_1}{B} < \frac{C}{B}$. This gives $B < C^3$. Similarly, one also has $A < C^3$.

\section{Proof of Theorem \ref{thm2}}
Without loss of generality, we can assume that $A$, $A+a_1$, $A+a_2$ are three consecutive divisors of $n$. From the proof of Theorem \ref{thm1}, we have $\frac{a_2 b_1 - b_2 a_1}{b_1 b_2} = \frac{a_2 - a_1}{B}$. Let $d = (a_2 - a_1) b_1 - (b_2 - b_1) a_1 \ge 1$. Then $B = \frac{b_2 b_1 (a_2 - a_1)}{d} = \frac{b_2 (a_1 (b_2 - b_1) + d)}{d} = \frac{1}{d} a_1 b_2 (b_2 - b_1) + b_2$. Suppose $a_2 - a_1 = \phi a_1$ and $b_2 - b_1 = \theta b_1$. Since $(a_2 - a_1) b_1 - (b_2 - b_1) a_1 > 0$, we have $\phi > \theta$.

\bigskip

Case 1: $a_2 = b_2$. Then from $A B = (A + a_2)(B - b_2)$, we have $a_2 = B - A > 0$. Now $A < A+a_1 < A+a_2 = B$ and $A + a_2 = B > B - b_1 > B - b_2 = A$. Since we assume that $A$, $A+a_1$, $A+a_2$ are three consecutive divisors of $n$, we must have $A+a_1 = B-b_1 =: M$ say. So $n = M^2 = (M + (a_2 - a_1))(M - (b_2 - b_1))$. This implies $(a_2 - a_1)(b_2 - b_1) = ((a_2 - a_1) - (b_2 - b_1))M$. As the left hand side is positive, the right hand side must be at least $M$. Therefore $C^2 > (a_2 - a_1)(b_2 - b_1) = ((a_2 - a_1) - (b_2 - b_1))M \ge M$ which implies $A, B < C^2 + C \le \frac{1}{4} C^3 - \frac{1}{2} C^2 + \frac{1}{4} C$ when $C \ge 7$.


\bigskip

Case 2: $|a_2 - b_2| = 1$. We will treat the case where $a_2 = b_2 + 1$. The other case $a_2 + 1 = b_2$ is similar. We have 
\begin{equation} \label{case2}
n = A B = (A + a_2)(B - (a_2 - 1))
\end{equation}
which implies $a_2 (B - A) + A = a_2 (a_2 - 1)$. So $a_2 | A$ and we write $A = a_2 A'$. Substituting this into (\ref{case2}), we have $A' B = (A' + 1)(B - (a_2 - 1))$ which implies $B - (a_2 - 1)A' = (a_2 - 1)$. So $a_2 - 1 | B$ and we write $B = (a_2 - 1) B'$. Substituting these new expressions for $A$ and $B$ into (\ref{case2}), we obtain $A' B' = (A' + 1)(B' - 1)$ which implies $B' = A' + 1$. Therefore
\[
A = a_2 A' \text{ and } B = (a_2 - 1) (A' + 1)
\]
and the three factorizations are
\[
n = [a_2 A'] [(a_2 - 1) (A' + 1)] = [a_2 A' + h] [(a_2 - 1) A' + k] = [a_2 A' + a_2] [(a_2 - 1) A']
\]
for some $0 < h < a_2$ and $0 < k < a_2-1$. Focusing on the last equality, we have, after some algebra,
\begin{equation} \label{case2.1}
(a_2 - h) (a_2 - 1 - k) A' = h k (A' + 1).
\end{equation}
Since $A'$ and $A' + 1$ are relatively prime, we must have $A' | hk$ and $A' + 1 | (a_2 - h)(a_2 - 1 - k)$.

\bigskip

Subcase 1: $h k = l A'$ and $(a_2 - h)(a_2 - 1 - k) = l (A'+1)$ for some integer $l \ge 2$. Then $l^2 A' (A' + 1) = h (a_2 - h) k (a_2 - 1 - k) \le (\frac{a_2}{2})^2 (\frac{a_2 - 1}{2})^2$. This implies $A' \le \frac{a_2 (a_2 - 1)}{4 l^2} \le \frac{a_2^2 - a_2}{16}$. So $A = a_2 A' \le \frac{a_2^3 - a_2^2}{16} \le \frac{1}{4} C^3 - \frac{1}{2} C^2 + \frac{1}{4} C$ for $C \ge 3$. Similarly one can show that $B \le \frac{1}{4} C^3 - \frac{1}{2} C^2 + \frac{1}{4} C$.

\bigskip

Subcase 2: $h k = A'$ and $(a_2 - h)(a_2 - 1 - k) = A' + 1$. Subtracting them, we have $a_2 (a_2 - 1) - h (a_2 - 1) - k a_2 = 1$ which implies $a_2 (a_2 - 1) - a_2 h - a_2 k = 1 - h$. So $a_2 | h - 1$ which forces $h = 1$ as $0 < h < a_2$. This implies $A' = k < a_2 - 1$ and hence $A \le a_2 (a_2 - 1) \le \frac{1}{4} C^3 - \frac{1}{2} C^2 + \frac{1}{4} C$ for $C \ge 6$. Similarly $B \le \frac{1}{4} C^3 - \frac{1}{2} C^2 + \frac{1}{4} C$.

\bigskip

Case 3: $a_2 - b_2 = 2$. Then
\begin{equation} \label{case3}
n = A B = (A + a_2) (B - (a_2 - 2))
\end{equation}
which implies $a_2 (B - A) - 2 A = a_2 (a_2 - 2)$. So $a_2 | 2A$ and write $2A = a_2 A'$. Multiplying (\ref{case3}) by two and substituting $2A = a_2 A'$, we have $A' B = (A' + 2)(B - (a_2 - 2))$ which implies $2B - (a_2 - 2)A' = 2 (a_2 - 2)$. We have $a_2 - 2 | 2B$ and write $2B = (a_2 - 2)B'$. Substituting these new expressions for $A$ and $B$ into (\ref{case3}), we have $A' B' = (A' + 2) (B' - 2)$ which implies $B' = A' + 2$. Therefore
\[
A = \frac{a_2 A'}{2} \text{ and } B = \frac{(a_2 - 2) (A' + 2)}{2}
\]
and the three factorizations are
\[
n = \Bigl[\frac{a_2 A'}{2} \Bigr] \Bigl[\frac{(a_2 - 2) (A' + 2)}{2} \Bigr] = \Bigl[\frac{a_2 A'}{2} + h \Bigr] \Bigl[\frac{(a_2 - 2) A'}{2} + k \Bigr] = \Bigl[\frac{a_2 (A' + 2)}{2} \Bigr] \Bigl[\frac{(a_2 - 2) A'}{2} \Bigr]
\]
for some $0 < h < a_2$ and $0 < k < a_2 - 2$. Focusing on the last equality, we have, after some algebra,
\begin{equation} \label{case3.1}
(a_2 - h)(a_2 - 2 - k) A' = h k (A' + 2).
\end{equation}

Subcase 1: If $A'$ is odd, then $A'$ and $A' + 2$ are relatively prime. Then we must have $A' | hk$ and $A' + 2 | (a_2 - h)(a_2 - 2 - k)$. Say $h k = l A'$ and $(a_2 - h)(a_2 - 2 - k) = l(A' + 2)$. So $l^2 A' (A' + 2) = h (a_2 - h) k (a_2 - 2 - k) \le (\frac{a_2}{2})^2 (\frac{a_2 - 2}{2})^2$ which implies $A' \le \frac{a_2 (a_2 - 2)}{4}$. Therefore $A = \frac{a_2 A'}{2} \le \frac{a_2^3 - 2a_2^2}{8} \le \frac{1}{4} C^3 - \frac{1}{2} C^2 + \frac{1}{4} C$ for $C \ge 2$. Similarly one also has $B \le \frac{1}{4} C^3 - \frac{1}{2} C^2 + \frac{1}{4} C$.

\bigskip

Subcase 2: If $A'$ is even, say $A' = 2 A''$. Then (\ref{case3.1}) becomes
\[
(a_2 - h)(a_2 - 2 - k) A'' = h k (A'' + 1).
\]
Since $A''$ and $A'' + 1$ are relatively prime, we must have $A'' | hk$ and $A''+2 | (a_2 - h)(a_2 - 2 - k)$. Say $hk = l A''$ and $(a_2 - h)(a_2 - 2 - k) = l (A'' + 1)$.

\bigskip

Subsubcase 1: If $l \ge 2$, then $l^2 A'' (A'' + 1) = h (a_2 - h) k (a_2 - 2 - k) \le (\frac{a_2}{2})^2 (\frac{a_2 - 2}{2})^2$ which implies $A'' \le \frac{a_2 (a_2 - 2)}{8}$. Hence $A' = 2A'' \le \frac{a_2 (a_2 - 2)}{4}$ and $A = \frac{a_2 A}{2} \le \frac{a_2^3 - 2 a_2^2}{8} \le \frac{1}{4} C^3 - \frac{1}{2} C^2 + \frac{1}{4} C$ for $C \ge 2$. Similarly $B \le \frac{1}{4} C^3 - \frac{1}{2} C^2 + \frac{1}{4} C$.

\bigskip

Subsubcase 2: If $l = 1$, then $h k = A''$ and $(a_2 - h)(a_2 - 2 - k) = A'' + 1$. Subtracting, we have $a_2(a_2 - 2) - (a_2 - 2)h - a_2 k = 1$. On one hand $a_2 (a_2 - 2) - a_2 h - a_2 k = 1 - 2h$ which implies $a_2 | 2h - 1$. Since $0 < h < a_2$, we must have $a_2 = 2h - 1$ which implies $h = \frac{a_2 + 1}{2}$. Substituting this into $a_2 (a_2 - 2) - a_2 h - a_2 k = 1 - 2h$, we obtain $k = \frac{a_2 - 3}{2}$. Therefore $A'' = (\frac{a_2 + 1}{2})(\frac{a_2 - 3}{2})$ and $A' = 2A'' = \frac{(a_2 + 1)(a_2 - 3)}{2}$. So $A = \frac{a_2 A'}{2} = \frac{a_2^3 - 2a_2^2 - 3a_2}{4} \le \frac{1}{4} C^3 - \frac{1}{2} C^2 + \frac{1}{4} C$ and $B = \frac{(a_2 - 2)(A' + 2)}{2} = \frac{(a_2 - 2)(a_2 - 1)^2}{4} = \frac{a_2^3 - 4a_2^2 + 5a_2 - 2}{4} \le \frac{1}{4} C^3 - \frac{1}{2} C^2 + \frac{1}{4} C$ for $C \ge 10$.

\bigskip

Case 4: $b_2 - a_2 = 2$. Then
\begin{equation} \label{case4}
n = A B = (A + (a_2-2)) (B - a_2)
\end{equation}
which implies $(a_2 - 2) B - a_2 A = a_2 (a_2 - 2)$. So $a_2 | 2B$ and write $2B = a_2 B'$. Multiplying (\ref{case4}) by two and substituting $2B = a_2 B'$, we have $A B' = (A + (a_2 - 2))(B' - 2)$ which implies $(a_2 - 2)B' - 2 A = 2 (a_2 - 2)$. We have $a_2 - 2 | 2A$ and write $2A = (a_2 - 2)A'$. Substituting these new expressions for $A$ and $B$ into (\ref{case4}), we have $A' B' = (A' + 2) (B' - 2)$ which implies $B' = A' + 2$. Therefore
\[
A = \frac{(a_2 - 2) A'}{2} \text{ and } B = \frac{a_2 (A' + 2)}{2}
\]
and the three factorizations are
\[
n = \Bigl[\frac{(a_2 - 2) A'}{2} \Bigr] \Bigl[\frac{a_2 (A' + 2)}{2} \Bigr] = \Bigl[\frac{(a_2 - 2) A'}{2} + h \Bigr] \Bigl[\frac{a_2 A'}{2} + k \Bigr] = \Bigl[\frac{(a_2 - 2) (A' + 2)}{2} \Bigr] \Bigl[\frac{(a_2 A'}{2} \Bigr]
\]
for some $0 < h < a_2 - 2$ and $0 < k < a_2$. Focusing on the last equality, we have, after some algebra,
\begin{equation} \label{case4.1}
(a_2 - 2 - h)(a_2 - k) A' = h k (A' + 2).
\end{equation}

Subcase 1: If $A'$ is odd, then $A'$ and $A' + 2$ are relatively prime. Then we must have $A' | hk$ and $A' + 2 | (a_2 - 2 - h)(a_2 - k)$. Say $h k = l A'$ and $(a_2 - 2 - h)(a_2 - k) = l(A' + 2)$. So $l^2 A' (A' + 2) = h (a_2 - 2 - h) k (a_2 - k) \le (\frac{a_2}{2})^2 (\frac{a_2 - 2}{2})^2$ which implies $A' \le \frac{a_2 (a_2 - 2)}{4}$. Therefore $A = \frac{(a_2 - 2) A'}{2} \le \frac{a_2^3 - 4a_2^2 + 4a_2}{8} \le \frac{1}{4} C^3 - \frac{1}{2} C^2 + \frac{1}{4} C$. Similarly one also has $B \le \frac{1}{4} C^3 - \frac{1}{2} C^2 + \frac{1}{4} C$ for $C \ge 4$.

\bigskip

Subcase 2: If $A'$ is even, say $A' = 2 A''$. Then (\ref{case4.1}) becomes
\[
(a_2 - 2 - h)(a_2 - k) A'' = h k (A'' + 1).
\]
Since $A''$ and $A'' + 1$ are relatively prime, we must have $A'' | hk$ and $A''+2 | (a_2 - 2 - h)(a_2 - k)$. Say $hk = l A''$ and $(a_2 - 2 - h)(a_2 - k) = l (A'' + 1)$.

\bigskip

Subsubcase 1: If $l \ge 2$, then $l^2 A'' (A'' + 1) = h (a_2 - 2 - h) k (a_2 - k) \le (\frac{a_2}{2})^2 (\frac{a_2 - 2}{2})^2$ which implies $A'' \le \frac{a_2 (a_2 - 2)}{8}$. Hence $A' = 2A'' \le \frac{a_2 (a_2 - 2)}{4}$ and $A = \frac{(a_2 - 2) A}{2} \le \frac{a_2^3 - 4 a_2^2 + 4 a_2}{8} \le \frac{1}{4} C^3 - \frac{1}{2} C^2 + \frac{1}{4} C$. Similarly $B \le \frac{1}{4} C^3 - \frac{1}{2} C^2 + \frac{1}{4} C$ for $C \ge 4$.

\bigskip

Subsubcase 2: If $l = 1$, then $h k = A''$ and $(a_2 - 2 - h)(a_2 - k) = A'' + 1$. Subtracting, we have $a_2(a_2 - 2) - (a_2 - 2) k - a_2 h = 1$. On one hand $a_2 (a_2 - 2) - a_2 h - a_2 k = 1 - 2k$ which implies $a_2 | 2k - 1$. Since $0 < k < a_2$, we must have $a_2 = 2k - 1$ which implies $k = \frac{a_2 + 1}{2}$. Substituting this into $a_2 (a_2 - 2) - a_2 h - a_2 k = 1 - 2k$, we obtain $h = \frac{a_2 - 3}{2}$. Therefore $A'' = (\frac{a_2 + 1}{2})(\frac{a_2 - 3}{2})$ and $A' = 2A'' = \frac{(a_2 + 1)(a_2 - 3)}{2}$. So $A = \frac{(a_2 - 2) A'}{2} = \frac{a_2^3 - 4a_2^2 + a_2 + 6}{4} \le \frac{1}{4} C^3 - \frac{1}{2} C^2 + \frac{1}{4} C$ for $C \ge 7$, and $B = \frac{a_2 (A' + 2)}{2} = \frac{a_2 (a_2 - 1)^2}{4} = \frac{a_2^3 - 2a_2^2 + a_2}{4} \le \frac{1}{4} C^3 - \frac{1}{2} C^2 + \frac{1}{4} C$.

\bigskip

Case 5: $|a_2 - b_2| \ge 3$. If $a_2 < b_2 \le C$, then $b_2 = (1 + \theta) b_1 \le C$ and $a_2 = (1 + \phi) a_1 \le C - 3$ Thus $a_1 \le \frac{C - 3}{1 + \phi} < \frac{C - 3}{1 + \theta}$, $b_1 \le \frac{C}{1 + \theta}$ and $b_2 - b_1 = \theta b_1 \le \frac{\theta C}{1 + \theta}$. Therefore, as $d \ge 1$,
\[
B = \frac{1}{d} a_1 b_2 (b_2 - b_1) + b_2 < \frac{C - 3}{1 + \theta} C \frac{\theta C}{1 + \theta} + C = \frac{\theta}{(1 + \theta)^2} C^2 (C - 3) + C.
\]
It remains to note that $\frac{\theta}{(1 + \theta)^2}$ has maximum value $\frac{1}{4}$ when $\theta = 1$. This gives $B < \frac{1}{4} C^2 (C - 3) + C \le \frac{1}{4} C^3 - \frac{1}{2} C^2 + \frac{1}{4} C$ for $C \ge 3$. When $b_2 < a_2 \le C$, one obtains $B < \frac{1}{4} C (C - 3)^2 + C \le \frac{1}{4} C^3 - \frac{1}{2} C^2 + \frac{1}{4} C$ via a similar argument. Analogously we have the same bound for $A$ which gives Theorem \ref{thm2}.

\section{Proof of Corollary \ref{cor1}}
Suppose $(x_1, y_1)$, $(x_2, y_2)$ and $(x_3, y_3)$ are three integer lattice points on the hyperbola $x y = n$ with $x_1 < x_2 < x_3$. Then $n = A B = (A + a_1)(B - b_1) = (A + a_2)(B - b_2)$ where $A = x_1$, $B = y_1$, $a_1 = x_2 - x_1$, $b_1 = y_1 - y_2$, $a_2 = x_3 - x_1$, $b_2 = y_1 - y_3$. Let $C = \max(a_2, b_2)$. Then Theorem \ref{thm2} tells us that $n = A B \le (\frac{1}{4} C^3 - \frac{1}{2} C^2 + \frac{1}{4} C)^2 < \frac{1}{16} (C - 0.5)^6$ which implies $2^{2/3} n^{1/6} + 0.5 < C$. We have the first half of Corollary \ref{cor1}. The second half of Corollary \ref{cor1} follows from the example in (\ref{eq1}). Recall
\[
n = N (N+1)^2 (N+2) (2N+1) (2N+3)
\]
where
\[
n = [(2N+1) N (N+2)] \cdot [(2N+3) (N+1)^2] =: x_1 \cdot y_1,
\]
\[
n = [(2N+1) N (N+2) + N] \cdot [(2N+3) (N+1)^2 - (N+1)] =: x_2 \cdot y_2,
\]
and
\[
n = [(2N+1) N (N+2) + (2N+1)] \cdot [(2N+3) (N+1)^2 - (2N+3)] =: x_3 \cdot y_3.
\]
Then $\max(x_3 - x_1, y_1 - y_3) = 2N+3$. One can see that $n/4 = N(N+1/2)(N+1)^2(N+3/2)(N+2) \ge (N + 0.9)^6$ when $n$ is sufficiently large. This implies
$n^{1/6} / 2^{1/3} > N + 0.9$ or $2^{2/3} n^{1/6} > 2N + 1.8$. Hence $2^{2/3} n^{1/6} + 1.2 > 2N+3$.

Department of Arts and Sciences \\
Victory University \\
255 N. Highland St., \\
Memphis, TN 38111 \\
U.S.A. \\
thchan@victory.edu

\end{document}